\documentclass[12pt]{article}
\usepackage{amsmath}
\usepackage{amssymb}
\usepackage{stmaryrd}
\usepackage{amsthm}

\usepackage[dvipsnames]{pstricks} 
\usepackage{multido}

\usepackage{graphicx}
\usepackage{xcolor}

\newcommand{\showgrid}{}

\newcommand{\gridon}{\renewcommand{\showgrid}{\psset{subgriddiv=1,griddots=10,gridlabels=6pt}\psgrid}}

\gridon 

\newgray{hellgrau}{.89}

\usepackage[applemac]{inputenc}

\usepackage{algorithmic}

\newif\ifenglish

\englishtrue

\newif\ifvariant
\variantfalse

\def\bit{\begin{itemize}}
\def\eit{\end{itemize}}
\def\beq{\begin{equation}}
\def\eeq{\end{equation}}

\englishtrue

\def\of#1{\left(#1\right)} 
\def\pas#1{\left(#1\right)} 
\def\setof#1{\left\{#1\right\}}

\def\abs#1{{\left\lvert{#1}\right\rvert}}

\def\strich{^\prime}

\def\defeq{\stackrel{\text{\tiny def}}{=}}

\def\N{{\mathbb N}}

\def\R{{\mathbb R}}

\def\0{{\mathbf 0}} 
\def\1{{\mathbf 1}} 

\def\CT1{CT1}
\def\xxxCT1{CT1}

\newtheoremstyle{excstyle}
  {1em}
  {2pt}
  {\sffamily\footnotesize\slshape}
  {0pt}
  {\sffamily\footnotesize\bfseries}
  {:}
  { }
  {}
\theoremstyle{excstyle}

\def\figref#1{\ifenglish Figure\else Abbildung\fi~\ref{#1}}

\ifenglish\else\fi

\def\figref#1{Fi\-gu\-re~\ref{#1}}
\def\defeq{:=}

\newgray{mfgray85}{0.85}

\newgray{mfgray70}{0.70}

\newgray{mfgray60}{0.60}

\title{A simple explanation for the ``shuffling phenomenon'' for lozenge tilings of dented hexagons}

\author{Markus Fulmek\thanks{
Research supported by the National Research Network ``Analytic
Combinatorics and Probabilistic Number Theory'', funded by the
Austrian Science Foundation. 
}\\
\small Fakult\"at f\"ur Mathematik \\
\small Oskar-Morgenstern-Platz 1, A-1090 Wien, Austria \\
\small \tt Markus.Fulmek@Univie.Ac.At
}

\date{2019--12--09}

\def\secA{\subsection*}

\def\EM#1{{\em #1\/}}
\begin{document}
\bibliographystyle{plain}


\long\def\insertfigure#1#2{%
\begin{figure}%
\begin{center}%
\input graphics/#1%
\end{center}%
\caption{#2}%
\label{fig:#1}%
\end{figure}%
}

\long\def\insertcommentedfigure#1#2#3{%
\begin{figure}%
\begin{center}%
\input graphics/#1%
\end{center}%
{#3}
\caption{#2}%
\label{fig:#1}%
\end{figure}%
}

\parindent0pt
\parskip1em

\def\lo{left--tilted}
\def\ro{right--tilted}
\def\vo{vertical}
\def\uo{upwards--pointing}
\def\do{downwards--pointing}
\def\gtp{GT\!P}

\maketitle

\begin{abstract}
In a recent paper, Lai and Rohatgi proved a ``shuffling theorem'' for lozenge tilings
of a hexagon with ``dents'' (i.e., missing triangles). Here, we shall point out that
this follows immediately from the enumeration of Gelfand--Tsetlin patterns with
given bottom row. This observation is also contained in a recent preprint of Byun.
\end{abstract}

\secA{Introduction}
In \cite[Theorems 2.1 and 2.3]{Lai:2019:ASTFLTODH}, Lai and Rohatgi observed that
there is a nice factorization of 
the quotient of enumerations of lozenge tilings for two ``hexagons with dents'',
where the hexagons only differ (in some sense to be explained below) by a ``shuffling of the dents'',
and asked for a
combinatorial explanation. 
  
The purpose of this note is to show that this observation follows from
the formula for the enumeration of Gelfand--Tsetlin patterns with prescribed bottom row.
I am grateful to Mihai Ciucu for pointing out to me that this observation is also
contained in a recent preprint of Byun \cite{Byun:2019:IISFATATAST}.

We will briefly  repeat the (well--known) background, mainly by
presenting pictures which illustrate the ideas.

\secA{Lozenge tilings of regions in the triangular lattice}

The \EM{triangular lattice} may be viewed as the \EM{tiling} of the plane
$\R^2$ by congruent equilateral triangles. A \EM{lozenge} (or \EM{rhombus}) is the
union of two such triangles which are \EM{adjacent} (i.e., sharing an edge).
A \EM{region} in the triangular lattice is simply a finite subset of triangles,
and a \EM{lozenge tiling} of such region $R$ is a partition of $R$ whose
blocks are lozenges (i.e., a family of disjoint lozenges whose union equals $R$).
The meaning of this dry set--theoretic definition becomes clear when looking
at some pictures; see \figref{fig:hexagonal-region}.

As usual in combinatorics, the first question is the \EM{enumeration} of such lozenge
tilings: It turns out that there is a wealth of interesting formulas for
regions of ``hexagonal shape''; see \figref{fig:hexagonal-region} again.

When visualizing the triangular lattice in the obvious way (see \figref{fig:hexagonal-region}),
it is immediately clear that there are precisely two possible orientations of the
equilateral triangles, namely \EM{\uo} and \EM{\do},
and that every lozenge consists of an \uo\ and a \do\
triangle, and that there are precisely three possible orientations of such lozenges,
namely \EM{\vo}, \EM{\lo} and \EM{\ro}.

So it is a necessary condition for the existence of lozenge tilings for some region
that the number of \uo\ triangles equals the number of \do\ triangles in this region
(whence the hexagon shown in \figref{fig:hexagonal-region} has \EM{no} lozenge tiling).

When considering some 
hexagonal shape with a \EM{horizontal diagonal} (i.e., a diagonal parallel to the
base of the hexagon, see \figref{fig:hexagonal-region}),
we may consider the bisection of 
the hexagon in an \EM{upper half}
and a \EM{lower half}: The upper half consists of \EM{rows of adjacent triangles} where
\uo\ and \do\ triangles alternate, 
and the same is true
for the lower half. In the upper half, these rows of triangles start and end with an
\uo\ triangle.  In the lower half, these rows of triangles start and end with a
\do\ triangle.
 Clearly, the lower half appears as upper half when reflected at the horizontal diagonal,
 and vice versa: We shall call either of these halves a \EM{half--hexagon}.

\begin{figure}%
\begin{center}%
\psset{unit=0.65cm}
\begin{pspicture}(-4.0,-0.0)(8.0,9)
\psset{linewidth=0.007,linestyle=solid,linecolor=lightgray}
\psline(0.0,0.0)(-2.0,3.4641)
\psline(1.0,0.0)(-1.5,4.3301)
\psline(2.0,0.0)(-1.0,5.1962)
\psline(3.0,0.0)(-0.5,6.0622)
\psline(4.0,0.0)(0.0,6.9282)
\psline(4.5,0.866)(0.5,7.7942)
\psline(5.0,1.7321)(1.0,8.6603)
\psline(5.5,2.5981)(2.0,8.6603)
\psline(6.0,3.4641)(3.0,8.6603)
\psline(4.0,0.0)(0.0,0.0)
\psline(4.5,0.866)(-0.5,0.866)
\psline(5.0,1.7321)(-1.0,1.7321)
\psline(5.5,2.5981)(-1.5,2.5981)
\psline(6.0,3.4641)(-2.0,3.4641)
\psline(5.5,4.3301)(-1.5,4.3301)
\psline(5.0,5.1962)(-1.0,5.1962)
\psline(4.5,6.0622)(-0.5,6.0622)
\psline(4.0,6.9282)(0.0,6.9282)
\psline(3.5,7.7942)(0.5,7.7942)
\psline(3.0,8.6603)(1.0,8.6603)
\psline(6.0,3.4641)(4.0,0.0)
\psline(5.5,4.3301)(3.0,0.0)
\psline(5.0,5.1962)(2.0,0.0)
\psline(4.5,6.0622)(1.0,0.0)
\psline(4.0,6.9282)(0.0,0.0)
\psline(3.5,7.7942)(-0.5,0.866)
\psline(3.0,8.6603)(-1.0,1.7321)
\psline(2.0,8.6603)(-1.5,2.5981)
\psline(1.0,8.6603)(-2.0,3.4641)
\psset{linewidth=0.05,linestyle=solid,linecolor=gray,linearc=0.03}
\pspolygon(0.0,0.0)(4.0,0.0)(6.0,3.4641)(3.0,8.6603)(1.0,8.6603)(-2.0,3.4641)
\psline[linestyle=dashed](-2,3.4641)(6,3.4641)
\rput(2,-0.4){{\tiny base (length $4$)}}
\rput(2,3.1){{\tiny horizontal diagonal (length $8$)}}
\rput(2,5.5){{\tiny upper half ($6$ rows)}}
\rput(2,1.2){{\tiny lower half ($4$ rows)}}
\psset{linewidth=0.03,linecolor=black,fillstyle=solid}
\pspolygon[fillcolor=Apricot](1.0,6.9282)(2.0,6.9282)(2.5,7.7942)(1.5,7.7942)
\pspolygon[fillcolor=Tan](1.5,7.7942)(2.5,7.7942)(2.0,8.6603)(1.0,8.6603)
\pspolygon[fillcolor=Mahogany](1.5,7.7942)(1.0,8.6603)(0.5,7.7942)(1.0,6.9282)
\end{pspicture}%
\end{center}%
{{\small The picture shows a \EM{hexagonal region} in the
\EM{triangular lattice}, where we started to construct a lozenge tiling at the
upper left corner of the hexagon: The three lozenges represent the three possible
orientations (\vo, \lo\ and \ro).
The
\EM{horizontal diagonal} (shown as a dashed line) bisects the hexagon in an upper and lower
\EM{half}. The upper half
consists of $6$ \EM{rows of triangles}, all of which start and end with an \EM{\uo\ triangle}, and the lower half consists of $4$ \EM{rows of triangles}, all of which start and end with a
\EM{\do\ triangle}}}
\caption{Hexagonal region with partial lozenge tiling.}%
\label{fig:hexagonal-region}%
\end{figure}%

\secA{Row--wise construction of lozenge tilings}
We start the construction of a lozenge tiling for some hexagon 
in the top row of the upper half: If such row contains $n$ \uo\ triangles, then there
are precisely $n-1$ \do\ triangles; and covering the latter by lozenges (which are
necessarily 
\lo\ or \ro) leaves precisely one
remaining \uo\ triangle, which then clearly \EM{must} be covered by a \vo\ lozenge
(see \figref{fig:rowwise}).

If we number the \uo\ triangles in rows from the left, then we may encode
this partial lozenge tiling by the number of this \vo\ lozenge. This lozenge
\EM{protrudes} to the \EM{next} row where it occupies one \do\ triangle, while
the remaining \do\ triangles there 
must be covered by \lo\ or \ro\ lozenges, 
leaving precisely two \uo\ triangle which \EM{must} be covered
by \vo\ lozenges (see \figref{fig:rowwise} again).

\begin{figure}%
\begin{center}%
\psset{unit=0.65cm}
\begin{pspicture}(-1.0,-4.0)(16.0,1.0)
\psset{linewidth=0.01,linecolor=lightgray,fillstyle=none}
\psline(0.5,0.0)(6.5,0.0)
\psline(1.0,0.866)(6.0,0.866)
\psline(0.5,0.0)(1.0,0.866)(1.5,0.0)
\psline(1.5,0.0)(2.0,0.866)(2.5,0.0)
\psline(2.5,0.0)(3.0,0.866)(3.5,0.0)
\psline(3.5,0.0)(4.0,0.866)(4.5,0.0)
\psline(4.5,0.0)(5.0,0.866)(5.5,0.0)
\psline(5.5,0.0)(6.0,0.866)(6.5,0.0)
\rput(1.0,0.2887){{\lightgray\tiny 1}}
\rput(2.0,0.2887){{\lightgray\tiny 2}}
\rput(3.0,0.2887){{\lightgray\tiny 3}}
\rput(4.0,0.2887){{\lightgray\tiny 4}}
\rput(5.0,0.2887){{\lightgray\tiny 5}}
\rput(6.0,0.2887){{\lightgray\tiny 6}}
\pspolygon[fillstyle=solid,fillcolor=gray](3.5,0.0)(4.0,0.866)(4.5,0.0)
\rput(4.0,0.2887){{\tiny\bf 4}}
\psline(8.5,0.0)(14.5,0.0)
\psline(9.0,0.866)(14.0,0.866)
\psline(8.5,0.0)(9.0,0.866)(9.5,0.0)
\psline(9.5,0.0)(10.0,0.866)(10.5,0.0)
\psline(10.5,0.0)(11.0,0.866)(11.5,0.0)
\psline(11.5,0.0)(12.0,0.866)(12.5,0.0)
\psline(12.5,0.0)(13.0,0.866)(13.5,0.0)
\psline(13.5,0.0)(14.0,0.866)(14.5,0.0)
\pspolygon[fillstyle=solid,fillcolor=Mahogany](11.5,0.0)(12.0,0.866)(12.5,0.0)(12.0,-0.866)
\psset{fillstyle=solid,linecolor=black,fillcolor=Apricot}
\pspolygon(8.5,0.0)(9.5,0.0)(10.0,0.866)(9.0,0.866)
\pspolygon(9.5,0.0)(10.5,0.0)(11.0,0.866)(10.0,0.866)
\pspolygon(10.5,0.0)(11.5,0.0)(12.0,0.866)(11.0,0.866)
\psset{fillstyle=solid,fillcolor=Tan}
\pspolygon(12.5,0.0)(13.5,0.0)(13.0,0.866)(12.0,0.866)
\pspolygon(13.5,0.0)(14.5,0.0)(14.0,0.866)(13.0,0.866)
\psset{fillstyle=none,linecolor=lightgray}
\psline(0.0,-3.0)(7.0,-3.0)
\psline(0.5,-2.134)(6.5,-2.134)
\psline(0.0,-3.0)(0.5,-2.134)(1.0,-3.0)
\psline(1.0,-3.0)(1.5,-2.134)(2.0,-3.0)
\psline(2.0,-3.0)(2.5,-2.134)(3.0,-3.0)
\psline(3.0,-3.0)(3.5,-2.134)(4.0,-3.0)
\psline(4.0,-3.0)(4.5,-2.134)(5.0,-3.0)
\psline(5.0,-3.0)(5.5,-2.134)(6.0,-3.0)
\psline(6.0,-3.0)(6.5,-2.134)(7.0,-3.0)
\rput(0.5,-2.7113){{\lightgray\tiny 1}}
\rput(1.5,-2.7113){{\lightgray\tiny 2}}
\rput(2.5,-2.7113){{\lightgray\tiny 3}}
\rput(3.5,-2.7113){{\lightgray\tiny 4}}
\rput(4.5,-2.7113){{\lightgray\tiny 5}}
\rput(5.5,-2.7113){{\lightgray\tiny 6}}
\rput(6.5,-2.7113){{\lightgray\tiny 7}}
\pspolygon[fillstyle=solid,fillcolor=gray](1.0,-3.0)(1.5,-2.134)(2.0,-3.0)(1.5,-3.866)
\rput(1.5,-2.7113){{\tiny\bf 2}}
\pspolygon[fillstyle=solid,fillcolor=gray](5.0,-3.0)(5.5,-2.134)(6.0,-3.0)(5.5,-3.866)
\rput(5.5,-2.7113){{\tiny\bf 6}}
\psline(0.5,-2.134)(6.5,-2.134)
\psline(1.0,-1.2679)(6.0,-1.2679)
\psline(0.5,-2.134)(1.0,-1.2679)(1.5,-2.134)
\psline(1.5,-2.134)(2.0,-1.2679)(2.5,-2.134)
\psline(2.5,-2.134)(3.0,-1.2679)(3.5,-2.134)
\psline(3.5,-2.134)(4.0,-1.2679)(4.5,-2.134)
\psline(4.5,-2.134)(5.0,-1.2679)(5.5,-2.134)
\psline(5.5,-2.134)(6.0,-1.2679)(6.5,-2.134)
\rput(1.0,-1.8453){{\lightgray\tiny 1}}
\rput(2.0,-1.8453){{\lightgray\tiny 2}}
\rput(3.0,-1.8453){{\lightgray\tiny 3}}
\rput(4.0,-1.8453){{\lightgray\tiny 4}}
\rput(5.0,-1.8453){{\lightgray\tiny 5}}
\rput(6.0,-1.8453){{\lightgray\tiny 6}}
\pspolygon[fillstyle=solid,fillcolor=gray](3.5,-2.134)(4.0,-1.2679)(4.5,-2.134)(4.0,-3.0)
\rput(4.0,-1.8453){{\tiny\bf 4}}
\psline(8.0,-3.0)(15.0,-3.0)
\psline(8.5,-2.134)(14.5,-2.134)
\psline(8.0,-3.0)(8.5,-2.134)(9.0,-3.0)
\psline(9.0,-3.0)(9.5,-2.134)(10.0,-3.0)
\psline(10.0,-3.0)(10.5,-2.134)(11.0,-3.0)
\psline(11.0,-3.0)(11.5,-2.134)(12.0,-3.0)
\psline(12.0,-3.0)(12.5,-2.134)(13.0,-3.0)
\psline(13.0,-3.0)(13.5,-2.134)(14.0,-3.0)
\psline(14.0,-3.0)(14.5,-2.134)(15.0,-3.0)
\pspolygon[fillstyle=solid,fillcolor=Mahogany](9.0,-3.0)(9.5,-2.134)(10.0,-3.0)(9.5,-3.866)
\pspolygon[fillstyle=solid,fillcolor=Mahogany](13.0,-3.0)(13.5,-2.134)(14.0,-3.0)(13.5,-3.866)
\psset{fillstyle=solid,linecolor=black,fillcolor=Apricot}
\pspolygon(8.0,-3.0)(9.0,-3.0)(9.5,-2.134)(8.5,-2.134)
\pspolygon(12.0,-3.0)(13.0,-3.0)(13.5,-2.134)(12.5,-2.134)
\psset{fillstyle=solid,fillcolor=Tan}
\pspolygon(10.0,-3.0)(11.0,-3.0)(10.5,-2.134)(9.5,-2.134)
\pspolygon(11.0,-3.0)(12.0,-3.0)(11.5,-2.134)(10.5,-2.134)
\pspolygon(14.0,-3.0)(15.0,-3.0)(14.5,-2.134)(13.5,-2.134)
\psset{fillstyle=none,linecolor=lightgray}
\psline(8.5,-2.134)(14.5,-2.134)
\psline(9.0,-1.2679)(14.0,-1.2679)
\psline(8.5,-2.134)(9.0,-1.2679)(9.5,-2.134)
\psline(9.5,-2.134)(10.0,-1.2679)(10.5,-2.134)
\psline(10.5,-2.134)(11.0,-1.2679)(11.5,-2.134)
\psline(11.5,-2.134)(12.0,-1.2679)(12.5,-2.134)
\psline(12.5,-2.134)(13.0,-1.2679)(13.5,-2.134)
\psline(13.5,-2.134)(14.0,-1.2679)(14.5,-2.134)
\pspolygon[fillstyle=solid,fillcolor=Mahogany](11.5,-2.134)(12.0,-1.2679)(12.5,-2.134)(12.0,-3.0)
\psset{fillstyle=solid,linecolor=black,fillcolor=Apricot}
\pspolygon(8.5,-2.134)(9.5,-2.134)(10.0,-1.2679)(9.0,-1.2679)
\pspolygon(9.5,-2.134)(10.5,-2.134)(11.0,-1.2679)(10.0,-1.2679)
\pspolygon(10.5,-2.134)(11.5,-2.134)(12.0,-1.2679)(11.0,-1.2679)
\psset{fillstyle=solid,fillcolor=Tan}
\pspolygon(12.5,-2.134)(13.5,-2.134)(13.0,-1.2679)(12.0,-1.2679)
\pspolygon(13.5,-2.134)(14.5,-2.134)(14.0,-1.2679)(13.0,-1.2679)
\psset{fillstyle=none,linecolor=lightgray}
\end{pspicture}%
\end{center}%
{{\small
The pictures show the simple idea of constructing
a lozenge tiling of some hexagon ``row by row''. Assume the top row of some
hexagon has length $6$. We may choose \uo\ triangle $4$ to be the
one which should be covered by a \vo\ lozenge. With this choice, the tiling of the remaining
top row is uniquely determined: The upper left pictures shows this top row
with the numbering of the \uo\ triangles, and the upper right picture shows
the (unique) lozenge tiling determined by the position ($4$) of the
\vo\ lozenge, which protrudes to the next row. 

In this next row,
we have to choose \EM{two} \uo\ pointing triangles which are covered by \vo\
lozenges, such that the protruding lozenge from the row above is located \EM{between} them:
The lower pictures show a possible choice of positions for such \uo\ triangles
($2$ and $6$) and the lozenge tiling uniquely determined by this choice.
Note that we may encode this partial lozenge tiling 
by the array of numbers $\pas{\pas{4},\pas{2,6}}$.
}}
\caption{Row--wise construction of lozenge tilings.}%
\label{fig:rowwise}%
\end{figure}%

\secA{Gelfand--Tsetlin patterns}
Clearly 
we can continue this process of ``row by row'' construction
of partial lozenge tilings as long as we are in the upper half of the hexagon, and
the partial tiling obtained by this construction is uniquely encoded by a triangular array of natural
numbers
$$
\begin{matrix}
 &  &  &  & u_{1,1} &  &  &  & \\
 &  &  & u_{2,1} & & u_{2,2} &  &  & \\
 &  & u_{3,1} & & u_{3,2} & & u_{3,3} &  & \\
 & u_{4,1} & & u_{4,2} & & u_{4,3} & & u_{4,4} & \\
\dots & & \dots & & \dots & & \dots & & \dots\\
\end{matrix}
$$
where $u_{i,j}$ indicates the \EM{number} of the $j$--th \uo\ triangle (counted from the left)
covered by a \vo\ lozenge in the $i$--th
row of triangles (counted from the top row).
The entries 
in row $i-1$ are in the following sense
``interlaced'' with the entries 
in row $i$:
\begin{equation}
\label{eq:interlacing}
u_{i,1}\leq u_{i-1,1}<u_{i,2}\leq u_{i-1,2} < u_{i,3}\leq\cdots u_{i-1,i-1}<u_{i,i}.
\end{equation}
Such arrays of natural numbers are called \EM{Gelfand--Tsetlin patterns}; see, for instance,
\cite{CohnLarsenPropp:1998:TSOATBPP}, where these patterns are considered in reverse order
(i.e., from bottom to top row) and called \EM{semi--strict Gelfand patterns}.
We shall also use the (space saving) one--line notation
$$\pas{
\pas{u_{1,1}},
\pas{u_{2,1},u_{2,2}},
\pas{u_{3,1},u_{3,2},u_{3,3}},
\pas{u_{4,1},u_{4,2},u_{4,3},u_{4,4}},
\dots,
}
$$
for Gelfand--Tsetlin patterns in the following.

Write $U_n=\setof{u_{n,1}<u_{n,2}<\cdots u_{n,n}}$
for the (ordered) set of the numbers in the $n$--the row, and denote the number of
Gelfand--Tsetlin patterns with bottom row $U_n$ as $\gtp\of{U_n}$.
From \eqref{eq:interlacing} we immediately obtain:
\begin{equation}
\label{eq:gtp-recursion}
\gtp\of{U_n}
=
\sum_{u_{n-1,1}=u_{n,1}}^{u_{n,2}-1}
\sum_{u_{n-1,2}=u_{n,2}}^{u_{n,3}-1}
\cdots
\sum_{u_{n-1,n-1}=u_{n,n-1}}^{u_{n,n}-1}
\gtp\of{U_{n-1}}.
\end{equation}
From \eqref{eq:gtp-recursion} and the initial condition $\gtp\of{U_1}\equiv1$ the following product formula is easily derived
(see
Cohn, Larsen and Propp 
\cite[Proposition 2.1]{CohnLarsenPropp:1998:TSOATBPP}):
\begin{equation}
\label{eq:gtp-formula}
\gtp\of{U_n} = \prod_{1\leq i<j\leq n}\frac{u_j-u_i}{j-i}.
\end{equation}

\secA{Putting together two half--hexagons}
\begin{figure}%
\begin{center}%
\input graphics/two-halves-equal-whole%
\end{center}%
{{\small
The left picture shows \EM{two} lozenge tilings
of a \EM{half--hexagon} with the \EM{same} positions of \vo\ lozenges
in the bottom row of triangles (one of the halves is reflected, so its bottom row now
apears as top row). These are encoded by Gelfand--Tsetlin patterns
$\pas{
\pas{4},
\pas{2,6},
\pas{1,4,7},
\pas{1,3,6,8}
}$ and 
$\pas{
\pas{3},
\pas{3,6},
\pas{2,5,7},
\pas{1,3,6,8}
}$ having the same bottom row.

Clearly, such halves can be combined 
to a lozenge tiling of the \EM{whole} hexagon, shown in the right picture.
}}
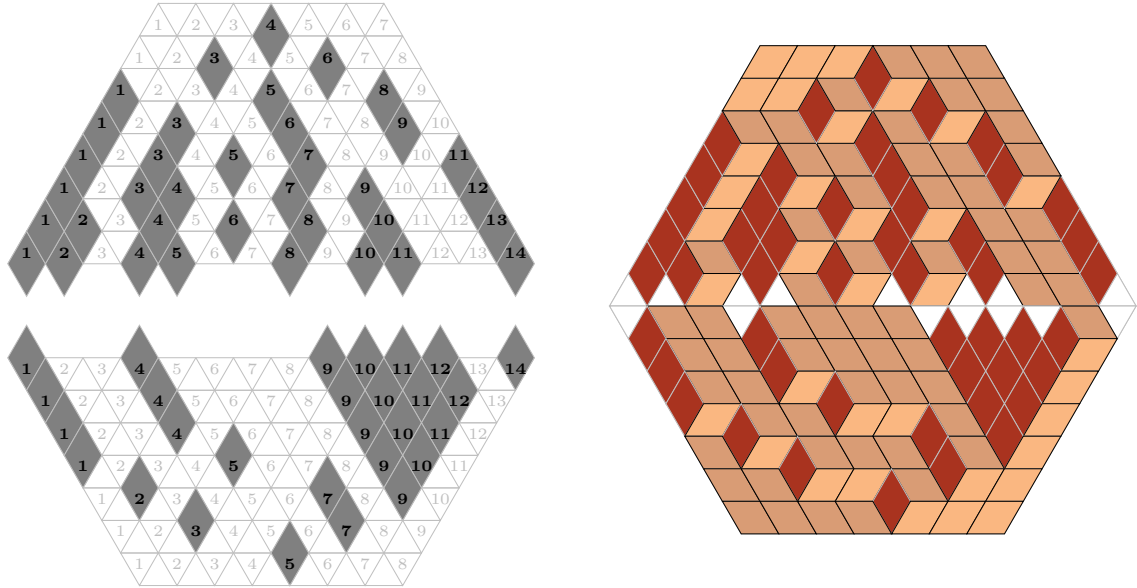
\caption{Putting together two half--hexagons.}%
\label{fig:two-halves-equal-whole}%
\end{figure}%

If we have two Gelfand--Tsetlin patterns \EM{with the same bottom row}, they
encode lozenge tilings for half--hexagons which we can \EM{put together}
to obtain a tiling of the ``whole'' hexagon (see \figref{fig:two-halves-equal-whole}).
Note that the hexagon thus obtained is \EM{not uniquely determined} by
the Gelfand--Tsetlin patterns: The only condition it must obviously fulfil is that
the \EM{length} of its horizontal diagonal must be greater or equal than the last
element in the bottom row of the Gelfand--Tsetlin patterns.

\secA{Dented hexagons}

\begin{figure}%
\begin{center}%
\input graphics/tri-lais-dented-hexagon%
\end{center}%
{{\small
Putting together the lozenge tilings encoded by the ``incompatible'' Gelfand--Tsetlin
patterns
{\tiny
$$\pas{
\pas{4},
\pas{3,6},
\pas{1,5,8},
\pas{1,3,6,9},
\pas{1,3,5,7,11},
\pas{1,3,4,7,9,12},
\pas{1,2,4,6,8,10,13},
\pas{1,2,4,5,8,10,11,14}
}
$$
}
and
{\tiny
$$\pas{
\pas{5},
\pas{3,7},
\pas{2,7,9},
\pas{1,5,9,10},
\pas{1,4,9,10,11},
\pas{1,4,9,10,11,12},
\pas{1,4,9,10,11,12,14}
}
$$}
(with \EM{different} bottom rows) is possible if we \EM{omit} the lozenges
bisected by the horizontal diagonal; the result of this operation is a lozenge
tiling of the \EM{dented hexagon} shown in the right picture, where the ``dents'' (i.e.,
the missing triangles adjacent to the horizontal triangles) are coloured white.
(This is the example presented in \cite[Figure 2.1]{Lai:2019:ASTFLTODH}.)
}}
\caption{Lozenge tiling of a dented hexagon.}%
\label{fig:tri-lais-dented-hexagon}%
\end{figure}%

For two \EM{arbitrary} Gelfand--Tsetlin patterns (with possibly \EM{different} bottom rows),
we may also put together the lozenge tilings of corresponding half--hexagons (see
\figref{fig:tri-lais-dented-hexagon}): In general, this procedure will not give a 
hexagon, but if we simply \EM{omit} the protruding \vo\ lozenges in the bottom rows,
we obtain a lozenge tiling of a hexagon with missing triangles (adjacent to the
horizontal diagonal). We shall call such region a \EM{dented hexagon}; by construction
it is determined by the bottom rows of the two Gelfand--Tsetlin patterns (but not uniquely:
As before, the length of the horizontal diagonal must be greater or equal than the
maximal entry in both bottom rows).

%

Denote the bottom rows of the Gelfand--Tsetlin patterns of the upper and lower
half of the dented hexagon by $U=\setof{u_1<u_2<\dots <u_m}$ and
$L=\setof{d_1<d_2<\dots< d_n}$, respectively. Then clearly the
number of lozenge tilings of this dented hexagon equals
$\gtp\of{U}\cdot \gtp\of{L}$. 

\secA{Simple observations}
For two finite sets of natural numbers $A$ and $B$ define
$$
P\of{A,B}\defeq\prod_{\stackrel{\pas{a,b}\in A\times B}{a<b}}\pas{b-a}
$$
and observe that we may rewrite \eqref{eq:gtp-formula} as follows: 
$$
\gtp\of{U} = \frac{P\of{U,U}}{2!\cdot3!\cdots\pas{\abs U-1}!}.
$$

Set $C=U\cap L$, $\overline{U}\defeq U\setminus C$ and $\overline{L}\defeq L\setminus C$.
Clearly, we have
$$
P\of{U,U} =
P\of{\overline U, \overline U}\cdot
P\of{\overline U, C}\cdot
P\of{C, \overline U}\cdot
P\of{C, C},
$$
and the analogous statement holds true for $P\of{L,L}$.
Hence the product $\gtp\of U\cdot \gtp\of L$ equals
\begin{equation}
\label{eq:enum-dented}
\frac{
P\of{\overline U, \overline U}\cdot
P\of{\overline L, \overline L}\cdot
P\of{\overline U\cup\overline L,C}\cdot
P\of{C,\overline U\cup\overline L}\cdot
P\of{C,C}^2
}
{
\pas{2!\cdot3!\cdots\pas{\abs U-1}!}\cdot\pas{2!\cdot3!\cdots\pas{\abs L-1}!}
}.
\end{equation}

Now fix some arbitrary subset 
$$
S\subseteq\overline U\cup\overline L.
$$
By the ``shuffling'' of dents corresponding to such set $S$ we mean that all elements
in $S$ belonging to $U$ are moved to $L$ and vice versa; i.e.,
the sets $U\strich$
and $L\strich$ obtained from this ``shuffling'' are
\begin{align*}
U\strich &= \pas{U\setminus S}\cup\pas{L\cap S}, \\
L\strich &= \pas{L\setminus S}\cup\pas{U\cap S}.
\end{align*}
Clearly, such ``shuffling'' does not change the set $\overline U\cup\overline L$:
$$
\overline{U\strich}\cup\overline{L\strich} = \overline U\cup\overline L.
$$
So 
from \eqref{eq:enum-dented} we obtain
{\small
$$
\frac{\gtp\of{U}\cdot \gtp\of{L}}{\gtp\of{U\strich}\cdot \gtp\of{L\strich}}=
\frac{
P\of{\overline U, \overline U}\cdot
P\of{\overline L, \overline L}
}
{
P\of{\overline{U\strich}, \overline{U\strich}}\cdot
P\of{\overline{L\strich}, \overline{L\strich}}
}
\cdot
\frac
{
\pas{2!\cdots\pas{\abs{U\strich}-1}!}\cdot\pas{2!\cdots\pas{\abs{L\strich}-1}!}
}
{
\pas{2!\cdots\pas{\abs U-1}!}\cdot\pas{2!\cdots\pas{\abs L-1}!}
}.
$$
}

\secA{Lai and Rohatgi's``shuffling phenomenon''}
Now fix some arbitrary subset $V\subseteq C = U\cap\ L$ and consider a dented
hexagon where the following triangles adjacent to the horizontal diagonal
have been removed:
\bit
\item \uo\ triangles at positions in $U\setminus V$
\item and \do\ triangles at positions in $L\setminus V$.
\eit

\EM{Every} lozenge tiling of this dented hexagon must contain $\abs V$ vertical
lozenges beaded along the horizontal diagonal, at positions \EM{outside} the set
$F\defeq\pas{U\cup L}\setminus V$.

We may consider the \EM{restriction} that the positions of the $\abs V$ vertical lozenges
must be chosen from some finite subset $B\subset\N\setminus F$ of natural numbers
and immediately obtain the following formula for the number of such ``restricted''
lozenge tilings of a dented hexagon given by $U$, $L$ and $B$ (again, such
hexagon is not unique: The length of its horizontal diagonal must be greater or equal than
$\max\of{U\cup L\cup B}$):
\begin{multline}
\label{eq:restricted-enumeration}
\frac{
P\of{\overline U, \overline U}\cdot
P\of{\overline L, \overline L}
}
{
\pas{2!\cdots\pas{\abs U-1}!}\cdot\pas{2!\cdots\pas{\abs L-1}!}
}
\cdot\\
\sum_{\stackrel{V\strich\subseteq B}{\abs{V\strich}=\abs V}}
P\of{\overline U\cup\overline L,X}\cdot
P\of{X,\overline U\cup\overline L}\cdot
P\of{X,X}^2,
\end{multline}
where the sets $X$ appearing in the summands are defined as
$$X=\pas{\pas{U\cap L}\setminus V} \cup V\strich.$$

Since the ``shuffling'' corresponding to some set $S\subseteq \overline U\cup\overline L$
 does neither change the set $U\cap L$ nor the set
$\overline U\cup\overline L$, the sum in \eqref{eq:restricted-enumeration}
is invariant under such ``shuffling'':
This explains the observations formulated as Theorems 2.1 and 2.3 in
\cite{Lai:2019:ASTFLTODH}.

\bibliography{/Users/mfulmek/Work/TeX/database}

\end{document}